\documentclass[a4paper,12pt]{article}
\usepackage{amssymb,amsmath,graphics,pstricks}
\newenvironment{dwd}{\par\noindent{\bf Proof.}}{\par\rightline{$\blacksquare$}}

\newtheorem{theo}{Theorem}
\newtheorem{lema}{Lemma}
\newtheorem{prop}{Proposition}
\newtheorem{defi}{Definition}

\newtheorem{coro}{Corollary}

\def\be#1\ee{\begin{equation}#1\end{equation}}
\newcommand{\ba}{\begin{eqnarray} }
\newcommand{\ea}{\end{eqnarray} }
\def\bt#1\et{\begin{theo}#1\end{theo}}
\def\bl#1\el{\begin{lema}#1\end{lema}}
\def\bp#1\ep{\begin{prop}#1\end{prop}}
\def\bd#1\ed{\begin{defi}#1\end{defi}}
\def\ra{\rightarrow}

\def\E{\mathbf{E}}

\def\N{{\mathbb N}}

\def\ls{\leqslant}
\def\gs{\geqslant}
\def\for{\mbox{for}}
\begin{document}
\title {A note on the Menchov-Rademacher Inequality}
\author{Witold Bednorz\\ \emph{University of Warsaw}}
\date{}
\maketitle
\begin{abstract}
We improve constants in the Rademacher-Menchov inequality by showing that
$$
\E(\sup_{1\ls k\ls n}|\sum^k_{i=1} X_i|^2)\ls (a+b\log^2_2 n),
$$
for all orthogonal random variables
$X_1,...,X_n$ such that $\sum^n_{k=1}\E|X_k|^2=1$. \\

{\bf Mathematical discipline:} probability theory

\end{abstract}

\section{Introduction}

\noindent
We consider real or complex orthogonal random variables $X_1,...,X_n$, i.e.
$$
\E |X_i|^2<\infty,\;\;1\ls i\ls n\;\;\mbox{and}\;\;
\E(X_iX_j)=0,\;\;1\ls i,j\ls n.
$$
Let us denote $S_j:=X_1+...+X_j$ for $1\ls j\ls n$, and $S_0=0$.
Clearly
$$
\E|S_j-S_i|^2=\sum^j_{k=i}\E |X_k|^2,\;\;\for\;i\ls j.
$$
The best constant in the Menchov-Rademacher inequality is defined by
$$
D_n:=\sup\E\sup_{1\ls i\ls n}|S_i|^2,
$$
where the supremum is taken over all orthogonal systems $X_1,...,X_n$,
which satisfy $\sum^n_{k=1}\E|X_k|^2=1$. We define also
$$
C:=\limsup_{n\ra \infty}\frac{D_n}{\log^2_2 n}.
$$
Rademacher \cite{R} in 1922 and indepenedently Menchov \cite{M} in 1923
proved that there exists $K>0$
such that for $n\gs 2$
$$
D_n\ls K\log_2^2 n,\;\;\mbox{hence}\;\; C\ls K.
$$
By now there are several different proofs of the above inequality.
The traditional proof of Rademacher-Menchov inequality
uses the bisection method (see Doob \cite{D},
and Lo{\'e}v \cite{L}), which leads to
$$
D_n\ls (2+\log_2 n)^2,\;\;n\gs 2,\;\;\mbox{hence}\;\;C\ls 1.
$$
In 1970 Kounias \cite{K} used a trisection method to get a finer inequality
$$
D_n\ls(\frac{\log_2 n}{\log_2 3}+2)^2,\;\;n\gs 2,\;\;\mbox{hence}\;\;
C\ls (\frac{\log_2 2}{\log_2 3})^2.
$$
S. Chobayan, S.Levental and H. Salehi \cite{C-L-S} proved
the following result
\be\label{r1}
D_{2n}\ls \frac{4}{3}D_n\;\;\mbox{if}\;\;D_n\ls 3;\;\; D_{2n}\ls ((D_n-\frac{3}{4})^{1/2}+\frac{1}{2})^2
\ee
and as a consequence they got the estimate $D_n\ls \frac{1}{4}(3+\log^2_2 n)$, $C\ls \frac{1}{4}$. 
An example given in \cite{C-L-S} shows that $D\gs \frac{\log^2_2 n}{\pi^2 \log^2_2 e}$ and thus 
$C\gs 0,04868$.
The aim of this paper is to improve the bisection method
and together with (\ref{r1}) to obtain that $C<\frac{1}{9}$.

\section{Results}

\bt\label{t1}
For each $n,m\in \N$ and $l>2$ the following inequality holds
$$
\sqrt{D_{n(2m+l)}}\ls \sqrt{D_n}+\sqrt{\max\{D_m,2D_{l-1}\}}.
$$
If $l=2$ then even stronger inequality holds true
$$
\sqrt{D_{n(2m+l)}}\ls \sqrt{D_n}+\sqrt{D_m}.
$$
\et
\begin{dwd}
Let us denote $p:=2m+l$. The triangle inequality yields
$$
|S_i|\ls |S_i-S_{pj}|+|S_{pj}|.
$$
Consequently
$$
\max_{1\ls i\ls pn}|S_i|\ls \max_{1\ls i\ls pn}\min_{0\ls j \ls n}|S_i-S_{pj}|+\max_{0\ls j\ls n}|S_{pj}|.
$$
Thus
$$
\E\max_{1\ls i \ls pn}|S_i|^2\ls \E(\max_{1\ls i\ls pn}\min_{0\ls j \ls n}|S_i-S_{pj}|+\max_{0\ls j\ls n}|S_{pj}|)^2.
$$
The definition of $D_n$ together with the classical norm inequality implies
$$
\sqrt{D_{pn}}\ls \sqrt{D_n}+\sqrt{\E\max_i \min_{0\ls j \ls n}|S_i-S_{pj}|^2}
$$
It remains to show that
\ba
&& \E\max_{1\ls i \ls pn} \min_{0\ls j \ls n}|S_i-S_{pj}|^2\ls \max\{D_m,2D_{l-1}\},\;\;\mbox{if}\;\;l>2\nonumber\\
&& \E\max_{1\ls i\ls pn} \min_{0\ls j \ls n}|S_i-S_{pj}|^2\ls D_m\;\;\mbox{if}\;\;l=2.\nonumber
\ea
Let us denote
\ba
&&  A_j:=\max\{|S_i-S_{pj}|:\; pj\ls i\ls pj+m \},\nonumber\\
&&  B_j:=\max\{|S_{p(j+1)}-S_i|:\; pj+m+l\ls i\ls p(j+1) \}\nonumber\\
&&  C_j:=\max\{|S_i-S_{pj+m}|:\; pj+m<i<pj+m+l\}\nonumber\\
&&  D_j:=\max\{|S_{pj+m+l}-S_i|:\;pj+m<i<pj+m+l\}\nonumber
\ea
for each $j\in \{0,...,n-1\}$.
Each $0\ls i \ls dn$ can be written in the form $i=pj+r$, where
$j\in\{0,...,n-1\}$, $r\in\{1,2,...,p\}$.
If $r\ls m$, then
$$
|S_i-S_{pj}|^2\ls A^2_j.
$$
If $r\gs m+l$
$$
|S_{p(j+1)}-S_{i}|^2\ls B^2_j.
$$
The last case is when $i=pj+m+r$, $r\in\{1,...,l-1\}$. Let us denote
\ba
&&P_j:=S_{pj+m}-S_{pj},\;\;V_j:=S_{pj+m+r}-S_{pj+m},\;\;\nonumber\\
&&Q_j:=S_{p(j+1)}-S_{pj+m+l},\;\;W_j:=S_{pj+m+l}-S_{pj+m+r}.\nonumber
\ea
Clearly ($i=pj+m+r$, $r\in\{1,...,l-1\}$)
$$
\min\{|S_i-S_{pj}|^2,|S_{p(j+1)}-S_i|^2\}=\min\{|P_j+V_j|^2,|Q_j+W_j|^2\}.
$$
For all complex numbers $a,b,c,d$ there is
$$
\frac{1}{2}|a+b|^{2}\ls |a|^2+|b|^2,\;\;\frac{1}{2}|c+d|^2\ls |c|^2+|d|^2.
$$
Since
$$
\min\{|a+b|^2,|c+d|^2\}\ls \frac{1}{2}|a+b|^2+\frac{1}{2}|c+d|^2
$$
we obtain that
$$
\min\{|a+b|^2,|c+d|^2\}\ls |a|^2+|b|^2+|c|^2+|d|^2.
$$
Hence
$$
\min\{|S_i-S_{pj}|^2,|S_{p(j+1)}-S_i|^2\}
\ls |P_j|^2+|Q_j|^2+|V_j|^2+|W_j|^2.
$$
and consequently for each $pj<i\ls p(j+1)$ the following inequality holds
$$
\min\{|S_i-S_{pj}|^2,|S_{p(j+1)}-S_i|^2\}\ls A^2_j+B^2_j+C^2_j+D^2_j.
$$
In fact we have proved that
$$
\E\max_{1\ls i\ls pn}\min_{0\ls j \ls n}|S_i-S_{2(m+1)j}|^2\ls \E\sum^{n-1}_{j=0}(A^2_j+B^2_j+C^2_j+D^2_j).
$$
Let us observe that
\ba
&& \E A^2_j\ls D_m\sum^{m}_{k=1} \E |X_{pj+k}|^2,\;\; \E B^2_j\ls D_m\sum^{m}_{k=1 }\E|X_{pj+m+l+k}|^2,\nonumber\\
&& \E (C^2_j+D^2_j)\ls D_{l-1}(\E|X_{pj+m+1}|^2+\E|X_{pj+m+l}|^2+2\sum^{l-1}_{k=2}\E|X_{pj+m+k}|^2)\nonumber\;\;,
\ea
Notice that if $l=2$ then
$$
\E (C^2_j+D^2_j)\ls D_{1}(\E|X_{pj+m+1}|^2+\E|X_{pj+m+1}|^2)
$$
Hence, if $l>2$
$$
\E\max_{1\ls i\ls pn}\min_{0\ls j \ls n}|S_i-S_{2(m+1)j}|^2\ls \max\{D_m,2D_{l-1}\}
$$
and if $l=2$
$$
\E\max_{1\ls i\ls pn}\min_{0\ls j \ls n}|S_i-S_{2(m+1)j}|^2\ls D_m.
$$
It ends the proof.
\end{dwd}
\begin{coro}\label{c1}
For each $n\gs m$ the following inequality holds
$$
D_n\ls D_m(2+\frac{\log_2 n-\log_2 m}{\log_2 (2m+2)})^2.
$$
\end{coro}
\begin{dwd}
Taking $l=2$ in Theorem \ref{t1} we obtain
$$
D_{m(2m+2)^k}\ls D_m(k+1)^2.
$$
For  each $n\gs m$ there exists $k\gs 0$ such that $m(2m+2)^{k-1}<n\ls m(2m+2)^k$.
Hence
$$
k<1+\frac{\log_2 n-\log_2 m}{\log_2(2m+2)}.
$$
Consequently
$$
D_{n}\ls D_m(2+\frac{\log_2 n-\log_2 m}{\log_2 (2m+2)})^2,
$$
\end{dwd} 
The result implies
$$
C=\limsup_{n\ra \infty}\frac{D_n}{\log^2_2 n}\ls \frac{D_m}{\log^2_2 (2m+2)}.
$$
Putting $l>2$ in Theorem \ref{t1} and proceeding we prove in the same way as in Corollary \ref{c1})
we get the following result.
\begin{coro}\label{c2} For each $l>2$ and $n\gs m$ the inequality holds true
$$
C\ls \frac{\max\{D_{m},2D_{l-1}\}}{\log^2_2(2m+l)}.
$$
\end{coro}
Let us remind that $D_2=4/3$. Hence applying Corollary \ref{c1} with $m=2$ we get
$$
C\ls\frac{4}{3 \log^2_2 6}<\frac{1}{5}.
$$
Observe that due to (\ref{r1})
$$
D_{2}=\frac{4}{3},\;\;D_4\ls (\frac{4}{3})^2,\;\;D_{8}\ls (\frac{4}{3})^3,\;\;D_{16}\ls
(\frac{4}{3})^4
$$
and
$$
D_{32}\ls ((\frac{4}{3})^4-\frac{3}{4})^{1/2}+\frac{1}{2})^2,\;\;
D_{64}\ls (((D_{32}-\frac{3}{4})^{1/2}+\frac{1}{2})^2.
$$
Hence
$$
D_8\ls 2,3704\;\; D_{64}\ls 5,5741.
$$
Applying Corollary \ref{c2}  with $m=64$, $c=9$ we obtain
$$
C\ls 0,1107<1/9.
$$ 

\flushright{
Witold Bednorz\\
Institute of Mathematics\\
Warsaw University\\
Banacha 2, 02-097\\
POLAND\\
e-mial: wbednorz@mimuw.edu.pl
}

\end{document}